\documentclass[12pt]{amsart}

\usepackage{geometry}                % See geometry.pdf to learn the layout options. There are lots.
\usepackage{graphicx}
\usepackage{amssymb}

\newtheorem{theorem}{Theorem}[section]

\newtheorem{lemma}[theorem]{Lemma}
\newtheorem{claim}[theorem]{Claim}
\newtheorem{corollary}[theorem]{Corollary}

\theoremstyle{definition}

\newtheorem{problem}[theorem]{Problem}

\theoremstyle{remark}
\newtheorem{remark}[theorem]{Remark}

\numberwithin{equation}{section}

\newcommand{\cA}{\mathcal{A}}

\newcommand{\cB}{\mathcal{B}}

\newcommand{\D}{\mathbb{D}}

\newcommand{\C}{\mathbb{C}}
\newcommand{\T}{\mathbb{T}}

\newcommand{\N}{\mathbb{N}}

\newcommand{\R}{\mathbb{R}}

\newcommand{\Om}{\Omega }

\newcommand{\diam}{\operatorname{Diam}}

\newcommand{\ima}{\operatorname{Im}}

\newcommand{\defeq}{\mathrel{\mathop:}=}

\title{On a theorem of Landau and Toeplitz}
\author{Robert B. Burckel}
\address{Department of Mathematics, Cardwell Hall, Kansas State University,
Manhattan, KS 66506, USA}
\email{burckel@math.ksu.edu}

\author{Donald E. Marshall}
\address{Department of Mathematics, Box 354350
University of Washington 
Seattle, WA 98195-4350 USA}
\email{marshall@math.washington.edu}

\author{Pietro Poggi-Corradini}
\address{Department of Mathematics, Cardwell Hall, Kansas State University,
Manhattan, KS 66506, USA}
\email{pietro@math.ksu.edu}
\date{March 20, 2006}     
\begin{document}
\begin{abstract}
  The now canonical proof of Schwarz's Lemma appeared in a 1907
  paper of Carath\'eodory, who attributed it to Erhard
  Schmidt. Since then, Schwarz's Lemma has acquired considerable
  fame, with multiple extensions and generalizations. Much less known
  is that, in the same year 1907, Landau and Toeplitz obtained a
  similar result where the diameter of the image set takes over the
  role of the maximum modulus of the function. We give a streamlined
  proof of this result and also extend it to include bounds on the
  growth of the maximum modulus.
\end{abstract}

\maketitle
\baselineskip=18pt

\section{Schwarz's Lemma}
First, let us set the following standard notations: $\C$ denotes the
complex numbers, $\D\defeq\{z\in \C: |z|<1\}$ is the open unit disk, and
$\T\defeq\{z\in \C: |z|=1\}$ is the unit circle. Moreover, for $r>0$, we let
$r\D\defeq\{z\in\C: |z|<r\}$ and $r\T\defeq\{z\in \C: |z|=r\}$.

We begin by recalling the aforementioned
\begin{theorem}[Schwarz's Lemma]\label{thm:sl}
Suppose $f$ is analytic on the unit disk $\D$ and
\begin{equation}\label{eq:max}
\sup_{|z|<1}|f(z)-f(0)|\leq 1.
\end{equation}

Then, 
\begin{equation}\label{eq:der}
|f^\prime(0)|\leq 1
\end{equation}
and 
\begin{equation}\label{eq:growth}
|f(z)-f(0)|\leq |z|
\end{equation}
for every $z\in \D$.

Moreover, equality holds in {\rm (\ref{eq:der})} or in  
{\rm (\ref{eq:growth})} at some point in
$\D\setminus\{0\}$ if and only if $f(z)=a+cz$ for some constants
$a,c\in \C$ where $|c|=1$.
\end{theorem}
The standard way to prove Schwarz's Lemma is to factor 
$f(z)-f(0)=zg(z)$, for some analytic function $g$, then
apply the maximum modulus theorem to $g$ to deduce that
$\sup_{|z|<1}|g(z)|\leq 1$. This argument first appeared in a paper of Carath\`eodory \cite{caratheodory1907} where the idea is attributed to E.~Schmidt.  See Remmert, \cite{remmert1991} p.~272-273, and Lichtenstein, \cite{lichtenstein1919} footnote 427, for historical accounts.

\section{The theorem of Landau and Toeplitz}

In a 1907 paper, Landau and Toeplitz prove a similar result where
(\ref{eq:max}) is replaced by the diameter of the image set. 
(For a set $E\subset\C$ the diameter is $\diam E\defeq\sup_{z,w\in E}|z-w|$.) 

\begin{theorem}[Landau-Toeplitz \cite{landau-toeplitz1907}]\label{thm:lt}
Suppose $f$ is analytic on the unit disk $\D$ and $\diam
f(\D)\leq 2$. Then  
\begin{equation}\label{eq:lt}
|f^\prime(0)|\leq 1.
\end{equation}

Moreover, equality holds in {\rm (\ref{eq:lt})} if and only if
$f(z)=a+cz$ for some $a,c\in \C$ where $|c|=1$.
\end{theorem}
\begin{remark}
As we will see in the proof below 
the inequality (\ref{eq:lt}) is a simple consequence of (\ref{eq:der}). The main hurdle is proving the case of equality. Inequality (\ref{eq:lt}) appears in the classic book of P\'olya and Szeg\"o, p.~151 and p.~356 \cite{polya-szego1972}, where the paper of Landau and Toeplitz is mentioned. However, the case of equality is not discussed.
\end{remark}
\begin{remark}
Notice, for instance, that Theorem \ref{thm:lt} covers the case when $f(\D)$ is an equilateral triangle of side-length $2$, which is of course not contained in a disk of radius $1$.
\end{remark}
\begin{proof}
Decompose $f$
into its odd part and its even part: $f(z)=f_o(z)+f_e(z)$, where
$f_o(z)\defeq (f(z)-f(-z))/2$ and $f_e(z)\defeq (f(z)+f(-z))/2$. Then
$|f_o(z)|\leq \diam f(\D)/2\leq 1$, $f_o(0)=0$ and
$f_o^\prime(0)=f^\prime(0)$, so by Schwarz's Lemma (Theorem
\ref{thm:sl}) applied to $f_o$ we get 
\[|f^\prime(0)|\leq 1.  \] 
Theorem \ref{thm:lt} is a consequence of  the following claim.
\begin{claim}\label{cl:lt} 
If $|f^\prime(0)|=1$, then $f(z)\equiv f(0)+f^\prime(0)z$.  
\end{claim}
From the `equality' part of Schwarz's Lemma (Theorem \ref{thm:sl}), we find that
\begin{equation}\label{eq:fo}
f_o(z)=\frac{f(z)-f(-z)}{2}=f^\prime(0)z\qquad\forall z\in \D
\end{equation}
We use (\ref{eq:fo}) to show that $f$ is linear.

For $0\leq r <1$, 
let $D_r\defeq \diam f(r\D)$. First, we show that $D_r=\diam
f(r\T)$. Indeed, notice that by the Open Mapping Theorem
the set 
\[\{f(z)-f(w): |z|<r,\ |w|<r\}=\bigcup_{|w|<r} [f(r\D)-f(w)]
\]
is open. Therefore, no number
in this set has modulus $D_r$. However, there are points $z_0,w_0\in
r\overline{\D}$ with $|f(z_0)-f(w_0)|=D_r$. So one at least of them,
say $w_0$, must lie on $r\T$. But then since $f(r\D)-f(w_0)$ is open,
$z_0$ cannot belong to $r\D$. Thus both $z_0,w_0$ lie on $r\T$, which
proves that $D_r=\diam f(r\T)$.

Now we show that the diameter of the
image grows linearly, more precisely,  $D_r=2r$ for every $0\leq r<1$.
 
Since 
\[
h_u(z)\defeq \frac{f(z)-f(-uz)}{z}
\] 
is an analytic
function for $z\in \D$ whenever $u\in \T$ is fixed, its
maximum modulus on 
the disk $r\overline{\D}$ is either constant in $r$ or increasing in $r$. So the quantity 
\[
\frac{D_r}{r}=
\frac{\diam
  f(r\T)}{r}=\max_{|z|=r}\max_{|u|=1}\left|\frac{f(z)-f(-uz)}{z}\right| 
=\max_{|u|=1}\max_{|z|\leq r}|h_u(z)|
\]
is also either constant in $r$ or increasing in $r$.
But, on one hand,
\[
\frac{D_r}{r}\geq\sup_{|u|=1}|h_u(0)|=\sup_{|u|=1}|1+u||f^\prime(0)|=2.
\]
And, on the other hand, since $D_r/r\leq D_1/r\leq 2/r$,
\[
\lim_{r\uparrow 1}\frac{D_r}{r}\leq 2.
\]
So $D_r/r$ must be constant and $D_r=2r$ for every $0<r<1$.

Now, for every $|w|<1$, consider the function
\[
g_w(z)\defeq\frac{f(z)-f(-w)}{2f^\prime(0)}
\]
which is analytic for $z\in \D$.
Then, by (\ref{eq:fo}), $g_w(w)=w$. Also, if $0<|w|=r<1$,
\[
|g_w(w)|=r=\frac{D_r}{2}\geq \sup_{|z|<r} |g_w(z)| 
\]
i.e., $g_w$ fixes $w$ and preserves the disk $D(0,|w|)$ centered at $0$ and of
radius $|w|$. Using Lemma \ref{lem:fp} below, we get that $\ima
g_w^\prime(w)=0$. Therefore,  
\begin{equation}\label{eq:pos}
\ima\frac{f^\prime(w)}{2f^\prime(0)}= 0\qquad\forall w\in \D,
\end{equation}
whence, thanks to the Open Mapping Theorem,
$f^\prime(w)/(2f^\prime(0))$ is constant and equal to $1/2$. Thus,
\[
f(z)\equiv f(0)+f^\prime(0)z\qquad \forall z\in \D,
\]
which proves Claim \ref{cl:lt} and hence Theorem \ref{thm:lt}.
\end{proof}

We are left to show the following lemma.
\begin{lemma}\label{lem:fp}
Suppose $g$ is analytic in $\D$, $0<r<1$, $|w|=r$ and 
\[
w=g(w)\mbox{ and }r=\max_{|z|=r}|g(z)|.
\]

Then, $\ima g^\prime(w)=0$.
\end{lemma}
\begin{proof}
Actually, the stronger conclusion $g^\prime(w)\geq 0$ is geometrically
obvious because when $g^\prime(w)\neq 0$, the map $g$ is very close to
the rotation-dilation centered at $w$ given by $\zeta\mapsto
w+g^\prime(w)(\zeta-w)$. But since $g$ can't rotate 
points inside $D(0,|w|)$ to a point outside, the derivative must be
positive.

For the sake of rigor, we instead give a ``calculus'' proof of the weaker
statement, which has the advantage of being more historically accurate,
since it can be
perceived in the original paper of Landau and Toeplitz, and which they credit
to F. Hartogs. 

For $\theta\in \R$ introduce
\[
\phi(\theta)\defeq
|g(we^{i\theta})|^2=g(we^{i\theta})\overline{g(we^{i\theta})}.
\]
The function $g^\star(z)\defeq\overline{g(\bar{z})}$ is also analytic
in $\D$, and $\phi$ may be written
\[
\phi(\theta)=g(we^{i\theta})g^\star(\bar{w}e^{-i\theta}),
\]
enabling us to compute $\phi^\prime(\theta)$ via the product and chain
rules. We get routinely,
\[
\phi^\prime(\theta)=-2\ima \left[we^{i\theta}g^\prime(we^{i\theta})\overline{g(we^{i\theta})}\right]
\]
and setting $\theta=0$,
\[
\phi^\prime(0)=-2\ima
\left[wg^\prime(w)\overline{g(w)}\right]
=-2\ima \left[wg^\prime(w)\overline{w}\right]
=-2|w|^2\ima g^\prime(w).
\]
Since $\phi$ realizes its maximum over $\R$ at $\theta=0$, we have
$\phi^\prime(0)=0$, so the preceding equality proves Lemma \ref{lem:fp}. 
\end{proof}

\section{Some corollaries of the Landau-Toeplitz Theorem}

\begin{corollary}
Suppose $f$ is analytic on $\D$.
If $D_r=\diam f(r\D)$ is linear in $r$, then $f$ is linear.
\end{corollary}
\begin{proof}
By continuity, $D_r=cr$ for some $c>0$. Then $g\defeq 2f/c$ satisfies $\diam
g(\D)=2$ and by rotating and translating we get  $0<g^\prime(0)\leq
1$, as in the initial invocation of Schwarz's Lemma in the proof of
the Landau-Toeplitz Theorem \ref{thm:lt}. However, for $r$ small
$g(r\D)$ is almost round, hence  
\[
2r=D_r=\diam g(r\D)=2rg^\prime(0)+o(r)
\]
i.e., $g^\prime(0)=1$. Now apply Claim \ref{cl:lt}.
\end{proof}
Recall that for a simply-connected planar domain $\Omega$ (not $\C$)
the {\sf hyperbolic density} $\rho_\Omega$ is defined on $\Omega$ so that  
\[
\rho_\Omega(w)|dw|=\rho_\Omega(f(z))|f^\prime(z)||dz|=\rho_\D(z)|dz|\defeq
\frac{|dz|}{1-|z|^2} 
\]
for some, and hence for every, conformal map $f$ of $\D$ onto $\Omega$. By Schwarz's Lemma (Theorem \ref{thm:sl}), the following monotonicity holds
\[
\Omega \subset \tilde{\Omega}\qquad\Longrightarrow\qquad \rho_\Omega(z) \geq\rho_{\tilde{\Omega}}(z) \qquad\forall z\in\Omega.
\]
Also, one can check that if $\Omega$ is bounded, then $\rho_\Omega$ always attains its minimum.
\begin{corollary}\label{cor:hd}
Every simply-connected planar domain $\Omega$ with $\diam \Om\leq 2$ satisfies 
\begin{equation}\label{eq:rho}
\min_{w\in \Omega}\rho_\Omega(w)\geq 1.
\end{equation}
Moreover, equality holds in {\rm (\ref{eq:rho})} if and only if $\Omega$ is a disk of radius $2$.
\end{corollary}
\begin{proof}
Fix $w\in \Omega$. By the Riemann mapping theorem there is a one-to-one and analytic map $f$ of $\D$ onto $\Omega$ such that $f(0)=w$.
Then $\rho_\Omega(w)=1/|f^\prime(0)|$. Now apply the Landau-Toeplitz Theorem \ref{thm:lt}.
\end{proof}

\section{Growth bounds}
In view of the growth bound (\ref{eq:growth}) in Schwarz's Lemma, it is natural to ask whether a similar statement holds in the context of `diameter'.
We offer the following result.

\begin{theorem}\label{thm:don}
Suppose $f$ is analytic on the unit disk $\D$ and $\diam f(\D)\leq 2$.
Then for all $z\in \D$
\begin{equation}\label{eq:don}
|f(z)-f(0)| \leq |z|\frac{2}{1+\sqrt{1-|z|^2}}.
\end{equation}
Moreover, equality holds in {\rm (\ref{eq:don})} 
at some point in $\D\setminus\{0\}$  
if and only if $f$ is a linear
fractional transformation of the form
\begin{equation}\label{eq:lft}
f(z)=c\frac{z-b}{1-\overline{b}z}+a
 \end{equation}
for some constants $a\in \C$, $b\in \D\setminus\{0\}$ and $c\in \T$.
\end{theorem}
\begin{remark}
In Schwarz's Lemma, equality in (\ref{eq:growth}) at some point
in $\D\setminus\{0\}$ holds if and only if equality holds at every
point $z\in \D$. This is not true any more in Theorem
\ref{thm:don}. Namely, when
$f$ is the linear fractional transformation in (\ref{eq:lft}), then equality in (\ref{eq:don})
occurs only for $z\defeq 2b/(1+|b|^2)$.  
\end{remark}
\begin{remark}
Since the origin does not play a special role, we can write
(\ref{eq:don}) more symmetrically as follows:
\[
\frac{|f(z)-f(w)|}{\diam f(\D)}\leq \frac{|z-w|}{|1-\bar{w}z|+\sqrt{(1-|z|^2)(1-|w|^2)}}
 \qquad\forall z,w \in \D.
\]
This is done by applying Theorem \ref{thm:don} to $f$ precomposed with
a M\"obius transformation, and using a well-known identity for the
pseudo-hyperbolic metric.
\end{remark}

\begin{proof} Fix $d\in\D$ such that $f(d)\neq f(0)$. Set 
$$g= c_1 f\circ T +c_2$$
where $T$ is a linear fractional
transformation of $\D$ onto $\D$ such that
$T(x)=d$, $T(-x)=0$, for some $x>0$  and $c_1$, $c_2$ are constants
chosen so that  $g(x)=x$ and $g(-x)=-x$. By elementary algebra
$$T(z)=\frac{d}{|d|}\frac{z+x}{1+xz}$$
where $x\defeq |d|/(1+\sqrt{1-|d|^2})$, 
$$c_1\defeq \frac{2x}{f(d)-f(0)}\qquad \hbox{ and }\qquad
c_2\defeq -x\frac{f(d)+f(0)}{f(d)-f(0)}.$$
Then 
\begin{equation}\label{eq:diam}
\diam {g(\D)}=|c_1|\diam f(\D)\leq
\frac{4}{|f(d)-f(0)|}\frac{|d|}{(1+\sqrt{1-|d|^2})}.
\end{equation}
We now prove that $\diam g(\D) \ge 2$ with equality if and only
if $g(z)\equiv z$.

Set $h(z)\defeq (g(z)-g(-z))/2.$ Then $h(x)=x$
and $h(-x)=-x$. Note also that $h(0)=0$ so that $h(z)/z$ is analytic
in the disk 
and has value $1$ at $x$ and hence by the maximum principle
$\sup_\D |h(z)|= \sup_\D|h(z)/z|\ge 1$ with equality only if $h(z)=z$
for all $z\in \D$. Since, by definition of $h$, $\diam g(\D)\geq 2\sup_\D
|h|$, we see that $\diam g(\D)\geq 2$ and then (\ref{eq:diam})
gives (\ref{eq:don}) for $z=d$.

If equality holds in (\ref{eq:don}) at some point in
$\D\setminus\{0\}$, then that point is an eligible $d$ for the
preceding discussion, and (\ref{eq:diam}) shows that $\diam g(\D)\leq
2$, while we have already shown that $\diam g(\D)\geq 2$. Thus
$\diam g(\D) =2$. Hence $\sup_{z\in\D}|h(z)|=1$ and therefore $h(z)\equiv z$. 
Since $h$ is
the odd part of $g$,  we have
$g^\prime(0)=h^\prime(0)=1$. Thus, 
by the Landau-Toeplitz Theorem \ref{thm:lt} applied to $g$, we find that  $g(z)\equiv g(0)+z$ and thus
\[
f(z)=\frac{1}{c_1}T^{-1}(z)+f(T(0)).
\]
Moreover, equality at $z=d$ in (\ref{eq:don}) says that $|f(d)-f(0)|=2x$, hence $|c_1|=1$.  Since $T$ is a M\"{o}bius transformation of $\D$, namely of the form
\[
\eta\frac{z-\xi}{1-\overline{\xi}z}
\]
for some constants $\xi\in \D$ and $\eta\in \T$,
its inverse is also of this form. Therefore, we conclude that $f$ can be written as in (\ref{eq:lft}).

Finally, if $f$ is given by (\ref{eq:lft}), then $2b/(1+|b|^2)\in \D\setminus\{0\}$, and one checks that equality is attained in (\ref{eq:don}) when $z$ has this value and for no other value in $\D\setminus\{0\}$.
\end{proof}
\section{Higher derivatives}
We finish with a result, due to Kalle Poukka in 1907, which is to be compared with the usual Cauchy estimates that one gets from the maximum modulus.
Interestingly, Poukka seems to have been the first student of 
Ernst Lindel\"{o}f, who is often credited with having founded the
Finnish school of analysis.  
\begin{theorem}[Poukka \cite{poukka1907}]\label{thm:poukka}
Suppose $f$ is analytic on $\D$. Then for all positive integers $n$ we have
\begin{equation}\label{eq:pk}
\frac{|f^{(n)}(0)|}{n!}\leq \frac{1}{2}\diam f(\D).
\end{equation}
Moreover, equality holds in {\rm (\ref{eq:pk})} for some $n$ if and
only if $f(z)=f(0)+cz^n$ for some constant $c$ of modulus $\diam f(\D)/2$.
\end{theorem}
\begin{proof}
Write $c_k\defeq f^{(k)}(0)/k!$, so that $f(z)=\sum_{k=0}^\infty c_k
z^k$, for every $z\in \D$. Fix $n\in \N$. For every $z\in\D$ 
\begin{equation}\label{eq:coef}
h(z)\defeq f(z)-f(ze^{i\pi/n})=\sum_{k=1}^\infty c_k(1-e^{i\pi k/n})z^k.
\end{equation}
Fix $0<r<1$ and notice that, by absolute and uniform convergence,
\begin{equation}\label{eq:parseval}
\sum_{k=1}^\infty|c_k|^2|1-e^{i\pi k/n}|^2 r^{2k}=\int_0^{2\pi} |h(re^{i\theta})|^2\frac{d\theta}{2\pi}\leq (\diam f(\D))^2.
\end{equation}
Therefore
\[
|c_k(1-e^{i\pi k/n})|r^{k}\leq \diam f(\D)
\]
for every $0<r<1$ and every $k\in \N$.
In particular, letting $r$ tend to $1$ and then setting $k=n$, we get $2|c_n|\leq\diam f(\D)$, which is (\ref{eq:pk}). 

If equality holds here, then letting $r$ tend to $1$ in (\ref{eq:parseval}) we get that all coefficients $c_k(1-e^{i\pi k/n})$ in
(\ref{eq:coef}) for $k\neq n$ must be $0$. Hence, $c_k=0$ whenever $k$
is not a multiple of $n$. Thus, $f(z)=g(z^n)$ for some analytic
function $g$ on $\D$. Moreover, $g^\prime(0)=c_n$ and $\diam
g(\D)=\diam f(\D)$. So, by Theorem \ref{thm:lt}, $g(z)=cz$ for some
constant $c$ with $|c|=\diam g(\D)$, and the result follows.

\end{proof}

\section{Further problems}

Here we discuss a couple of problems that are related to these ``diameter" questions.

The first problem arises 
when trying to estimate the distance of $f$ from its linearization, $f(z)-(f(0)+f^\prime(0)z)$, to give a ``quantitative" version for the `equality' case in Schwarz's Lemma (Theorem \ref{thm:sl}).
This is done via the so-called Schur algorithm. As before, one considers the function
\[
g(z)\defeq \frac{f(z)-f(0)}{z}
\] 
which is analytic in $\D$, satisfies $g(0)=f^\prime(0)$ and which, by assumption (\ref{eq:max}) and the Maximum Modulus Theorem, has $\sup_\D|g|\leq 1$. Now let $a\defeq f^\prime(0)$ and post-compose $g$ with  a M\"{o}bius transformation of $\D$ which sends $a$ to $0$ to find that 
\[
\frac{g(z)-a}{1-\bar{a}g(z)}=zh(z)
\]
for some analytic function $h$ with $\sup_\D|h|\leq 1$.

Inserting the definition of $g$ in terms of $f$ and solving for $f$ shows that
\[
f(z)-f(0)-az=(1-|a|^2)\frac{z^2h(z)}{1+\bar{a}zh(z)}
\]
Thus, for every $0<r<1$,
\begin{equation}\label{eq:schur}
\max_{|z|<r}|f(z)-f(0)-f^\prime(0)z|\leq (1-|f^\prime(0)|^2)\frac{r^2}{1-|f^\prime(0)|r} 
\end{equation}
and `equality' holds for at least one such $r$ if and only if $h(z)\equiv a/|a|=f^\prime(0)/|f^\prime(0)|$, i.e., if and only if
\[
f(z)=z\frac{a}{|a|}\frac{z+|a|}{1+|a|z}+b
\]
for constants $a\in\overline{\D}$, $b\in \C$.
%
%Namely, suppose $f(z)=\sum_{n=0}^\infty a_nz^n$ is analytic in $\D$ and (\ref{eq:max}) holds. Then, by the same computation of the integral of the $|f(z)-f(0)|^2$ on the circle of radius $r$, $0<r<1$, as done in the proof of Theorem \ref{thm:poukka} above, we find that 
%\[
%\sum_{n=1}^\infty |c_n|^2\leq 1
%\]
%Thus, by Cauchy-Schwarz,
%\begin{eqnarray}
%|f(z)-(f(0)+f^\prime(0)z)|^2 & = & \left|\sum_{n=2}^\infty c_n z^n\right|^2\nonumber\\
%& \leq & \sum_{n=2}^\infty |c_n|^2 \sum_{n=2}^\infty |z|^{2n} \nonumber\\
%& \leq & (1-|f^\prime(0)|^2)\frac{|z|^4}{1-|z|^2}\label{eq:pom}
%\end{eqnarray}

In the context of this paper, when $f$ is analytic in $\D$ and $\diam f(\D)\leq 2$,
by the Landau-Toeplitz Theorem \ref{thm:lt} and a normal family argument we see that, for every $\epsilon>0$ and every $0<r<1$, there exists $\alpha>0$ such that: $|f^\prime(0)|\geq 1-\alpha$ implies 
\[
|f(z)-(f(0)+f^\prime(0)z)|\leq \epsilon\qquad \forall |z|\leq r.
\]
However, one could ask for an explicit bound as in (\ref{eq:schur}).
\begin{problem}\em
If $f$ is analytic in $\D$ and $\diam f(\D)\leq 2$, 
find an explicit (best?) function $\phi(r)$ for $0\leq r<1$ so that
\[
|f(z)-(f(0)+f^\prime(0)z)|\leq (1-|f^\prime(0)|)\phi(r)\qquad \forall |z|\leq r.
\]
\end{problem}

Another problem can be formulated in view of Corollary \ref{cor:hd}.
It is known, see the Corollary to Theorem~3 in \cite{minda-wright1982}, that if $\Omega$ is a bounded convex domain, then the minimum  
\begin{equation}\label{eq:min}
\Lambda(\Omega)\defeq \min_{w\in \Omega}\rho_\Omega(w)
\end{equation}
is attained at a unique point $\tau_\Omega$, which we can call the {\sf hyperbolic center} of $\Omega$. Also let us define the {\sf hyperbolic radius} of $\Omega$ to be
\[
R_h(\Omega)\defeq\sup_{w\in \Omega}|w-\tau_\Omega|.
\]

Now assume that
$\diam \Omega=2$. Then we know, by Corollary \ref{cor:hd}, that $\Lambda(\Omega)\geq 1$  with equality if and only if $\Omega$ is a disk of radius $1$. In particular, if $\Lambda(\Omega)=1$, then $R_h(\Omega)=1$.
\begin{problem}\em
Given $m>1$, find or estimate, in terms of $m-1$,
\[
\sup_{\Omega\in \cA_m} R_h(\Omega)
\]
where $\cA_m$ is the family of all convex domains $\Omega$ with $\diam \Omega=2$ and $\Lambda(\Omega)\leq m$.
\end{problem}

More generally, given an analytic function $f$ on $\D$ such that
$\diam f(\D)\leq 2$, define  
\[
M(f)\defeq\min_{w\in\D}\sup_{z\in\D}|f(z)-f(w)|
\]
and let $w_f$ be a point where $M(f)$ is attained. 
\begin{problem}\em
Fix $a<1$. Find or estimate, in terms of $1-a$,
\[
\sup_{f\in \cB_a}M(f)
\]
where $\cB_a$ is the family of all analytic functions $f$ on $\D$ with
$\diam f(\D)\leq 2$ and \[|f^\prime(w_f)|(1-|w_f|^2)\geq a.\]
\end{problem}

\def\cprime{$'$}

\end{document}